\documentclass{article}
\usepackage{verbatim, graphicx, caption, tikz, tikz-3dplot, amsmath, mathrsfs, array, tabu, siunitx, breqn, xcolor, hyperref, caption,imakeidx,pgfplots,enumitem,fancyhdr}
\pgfplotsset{compat=1.9}
\hypersetup{
    colorlinks,
    linkcolor={red!80!black},
    citecolor={blue!50!black},
    urlcolor={blue!80!black}
}
\usetikzlibrary{trees}
\usepackage[utf8]{inputenc}
\usepackage[linesnumbered,lined,ruled]{algorithm2e}
\usepackage[english]{babel}
\usepackage[]{amsthm} 
\usepackage[]{amssymb} 
\newcommand{\+}[1]{\ensuremath{\mathbf{#1}}}

\newcommand{\eps}{\varepsilon}
\newcommand{\del}{\delta}

\newcommand{\pde}{\partial}
\newcommand{\lam}{\lambda}
\newcommand{\aaa}{\noindent}

\newcommand*{\inv}[1]{{#1}^{\mkern-1.5mu\mathsf{-1}}}

\def\been{\begin{enumerate}}
\def\bean{\begin{enumerate}[label=(\alph*)]}
\def\enen{\end{enumerate}}
\def\beit{\begin{itemize}}
\def\enit{\end{itemize}}
\def\bbm{\begin{bmatrix}}
\def\ebm{\end{bmatrix}}
\def\bmm{\begin{matrix}}
\def\emm{\end{matrix}}
\def\bpm{\begin{pmatrix}}
\def\epm{\end{pmatrix}}
\newcommand{\bcc}{\begin{cases}}
\newcommand{\ecc}{\end{cases}}

\DeclareMathOperator{\diag}{diag}

\newcommand{\paren}[1]{\left(#1\right)}

\newcommand{\set}[1]{\left\{ #1 \right\}}

\newcommand{\abs}[1]{\left| #1 \right|}
\newcommand{\norm}[1]{\left\| #1 \right\|}
\newcommand{\bol}[1]{\textbf{#1}\index{#1}}

\usepackage{algorithmic}
\newcommand{\algorithmicinput}{\textbf{input}}
\newcommand{\INPUT}{\item[\algorithmicinput]}
\newcommand{\algorithmicoutput}{\textbf{output}}
\newcommand{\OUTPUT}{\item[\algorithmicoutput]}

\title{A Splitting Approach to Dynamic Mode Decomposition of Nonlinear Systems}

\author{Jovan \v Zigi\' c\\
\small{Dominican University New York, Orangeburg NY 10962, USA}}
\date{\today}

\begin{document}
\maketitle 


\begin{abstract}
Reduced-order models have long been used to understand the behavior of nonlinear partial differential equations (PDEs). Naturally, reduced-order modeling techniques come at the price of computational accuracy for a decrease in computation time. Optimization techniques are studied to improve either or both of these objectives and decrease the total computational cost of the problem. This paper focuses on the dynamic mode decomposition (DMD) applied to nonlinear PDEs with periodic boundary conditions. It provides a study of a newly proposed optimization framework for the DMD method called the Split DMD.
\end{abstract} 

\textbf{AMS subject classifications}: Primary 65K10, Secondary 65M22.

\section{Introduction}

The Navier-Stokes (NS) equations are the primary mathematical model for understanding the behavior of fluids. The existence and smoothness of the NS equations \cite{feff} is considered to be one of the most important open problems
in mathematics, and challenges in their numerical simulation is a barrier to understanding the physical phenomenon of turbulence. Due to the difficulty of studying this problem directly, problems in the form of nonlinear partial differential equations that exhibit similar properties to the NS equations are studied as preliminary steps towards building a wider understanding of the field. In fact, the use of the proper orthogonal decomposition method to build reduced-order models for fluids was accelerated with the publication of the first edition of the monograph ``Turbulence, Coherent Structures, Dynamical Systems and Symmetry" \cite{podbook}. 

Reduced-order models has since evolved as its own discipline with applications to control \cite{r1a}\cite{r1b}\cite{r1c}, optimization \cite{r2a}\cite{r2b}, and uncertainty quantification \cite{r3}. Naturally, reduced-order modeling techniques come at the price of computational accuracy for a decrease in computation time. Optimization techniques are studied to improve either or both of these objectives and decrease the total computational cost of the problem.
\section{Split Dynamic Mode Decomposition}

By definition of the term \cite{walter} and observation of simulated models governed by nonlinear PDEs, bifurcation present in system dynamics greatly alters system behavior over time. When approaching such a problem with DMD, this bifurcation effect leads to issues arising with selecting DMD modes that represent both {transient dynamics} and states of equilibrium that follow. This phenomena is a natural obstacle to both the DMD algorithm \cite{dmdbook} and the Levenberg-Marquardt algorithm \cite{lev44} \cite{marq63} crucial to the optimized DMD (OD) method \cite{ask17} \cite{optdmd}. 

In the DMD algorithm, the primary issue generated by effects of bifurcation is selecting DMD modes that accurately represent the dominant coherent structures, or pattern of behavior, of the dynamical system over the time interval. In a chaotic dynamical system, this pattern between the transient state and equilibrium state is not explicit. Thus, the linear problem that provides a solution for the low-rank Koopman operator is forced to compensate for the information provided throughout the entire interval, and consequently choose DMD modes that represent two distinct patterns. Thus, the fact that the resulting reduced-order system modeled by the DMD modes provides one pattern, or one set of coherent structures, to represent multiple distinct patterns is inherently flawed.

In the Levenberg-Marquardt algorithm (LMA), the primary issue generated by effects of transient states is the increased difficulty to solve the nonlinear least squares problem \cite{per18}, which directly contributes to computational complexity at the expense of time. Since approximating nonlinear PDEs presents nonconvex optimization problems through the LMA, globalization strategies such as trust region or line search methods must be employed to find stationary points. Generally, less iterations are required to solve the DMD problems when they contain simpler dynamics. However, when solving subproblems within each iteration becomes an increasingly difficult task, the cost of computational time becomes significant to the point of computational infeasibility.

This section proposes the \bol{split DMD} algorithm as a solution to the issue of selecting DMD or OD modes for nonlinear PDEs whose coherent structures vary over time. The split DMD algorithm uses a simple modification prior to the DMD routine to choose a domain for evaluation without changing the existing DMD method of computing system modes. The contribution to the existing DMD procedure comes from splitting the entire time interval into several subintervals before computing the low-rank Koopman modes. The \bol{split lines}, determining the boundaries of each subinterval, are selected by the \bol{n-split algorithm} described in the next section. The $n$-split algorithm produces $n$ number of subintervals by a recursive method. 

The goal of the algorithm is to separate the simulation interval into subintervals that have different ranges of data values, and are hypothesized to represent the time intervals with distinct system coherent structures. By heuristic reasoning, $n$ number of split lines are chosen as an initial guess for range comparison. These split lines are defined by the set of points in time:
\begin{align}\label{tsplit}
t_{split} = \set{t_0,t_1,t_2,\dotsc,t_n}
\end{align}
where a time interval $[0,T]$ implies $t_0 = 0$ and $t_n=T$,
so that the data (or system values) $Z(x,t)$ with respect to time $t$ and space $x$ is separated into subintervals:
\begin{equation}\label{eq:Zsplit}
 Z_{split} = \set{Z(x, [t_0, t_1]),Z(x, [t_1, t_2]),\dotsc,Z(x, [t_{n-1}, t_n])}
\end{equation}
By random selection, a line of space $x$ denoted $x_{test}$ is chosen to compare the system behavior across subintervals of time. Two tests are used to determine whether or not to merge two adjacent intervals. These involve either the difference in data value range or the length of the subinterval.

\section{Implementation of the Split DMD Algorithm}

For the first test, called the $\eps$ test, the necessary difference in range to retain a split line selection is data that exceeds a lower bound tolerance $\eps$, which is heuristically chosen as a fraction of the data range.  Furthermore, the minimum change $\eps$ in data is inspected at both upper and lower bounds of the subinterval data ranges, as illustrated in the following figures:
\begin{center}
\begin{minipage}{.5\textwidth}
\centering
\begin{tikzpicture}
 \draw[->] (4.0, -1.7) -- (8.0, -1.7) node[right] {$t$};
  \draw[->] (4, -1.7) -- (4, 1.3) node[above] {$z$};
\filldraw[black] (6,-1.7) circle (2pt) node[anchor=north] {$t_k$};
\filldraw[black] (5,-1.7) circle (2pt) node[anchor=north] {$t_{k-1}$};
\filldraw[black] (7,-1.7) circle (2pt) node[anchor=north] {$t_{k+1}$};
\draw (6.5,0.9) node[anchor=south] {$< \eps$};
\draw (5.5,-1.0) node[anchor=north] {$> \eps$};
\draw[gray, thick] (5,-1) -- (6,-1) -- (6,1) -- (5,1) -- (5,-1);
\draw[gray, thick] (6,-1.4) -- (7,-1.4) -- (7,0.9) -- (6,0.9) -- (6,-1.4);
\end{tikzpicture}
	\captionof{figure}{Retained Split Line $t_k$}
 	\label{fig:splitrange1}
\end{minipage}%
\begin{minipage}{.5\textwidth}
\centering
\begin{tikzpicture}
 \draw[->] (-2.0,-1.7) -- (2.0, -1.7) node[right] {$t$};
  \draw[->] (-2, -1.7) -- (-2, 1.3) node[above] {$z$};
\filldraw[black] (0,-1.7) circle (2pt) node[anchor=north] {$t_k$};
\filldraw[black] (-1,-1.7) circle (2pt) node[anchor=north] {$t_{k-1}$};
\filldraw[black] (1,-1.7) circle (2pt) node[anchor=north] {$t_{k+1}$};
\draw (0.5,0.9) node[anchor=south] {$< \eps$};
\draw (-0.5,-1.0) node[anchor=north] {$< \eps$};
\draw[gray, thick] (-1,-1) -- (0,-1) -- (0,1) -- (-1,1) -- (-1,-1);
\draw[gray, thick] (0,-1.1) -- (1,-1.1) -- (1,0.9) -- (0,0.9) -- (0,-1.1);
\end{tikzpicture}
	\captionof{figure}{Discarded Split Line $t_k$}
	 \label{fig:splitrange2}
\end{minipage}
\end{center}

\vspace{10pt}

In Figures \ref{fig:splitrange1} and \ref{fig:splitrange2}, the two-dimensional shapes represent the change in range of system values in time at a spatial point $x$.
For an example of the $\eps$ test, if the data range is $[0, 1]$, a chosen 10\% tolerance or $\eps = 0.1$ change in either bound between adjacent subintervals is the minimum amount required to retain the split line. Otherwise, if the adjacent subintervals for $k \in  \set{1,2,\dotsc,{n-1}}$ satisfy both of the following relationships:
\begin{equation}\label{eq:max}
  \abs{\max_{t \in [t_{k-1}, t_{k}] } Z(x_{test}, t) - \max_{t \in [t_{k}, t_{k+1}] }Z(x_{test}, t)} < \eps
\end{equation}
\begin{equation}\label{eq:min}
 \abs{\min_{t \in [t_{k-1}, t_{k}] } Z(x_{test}, t) - \min_{t \in [t_{k}, t_{k+1}] }Z(x_{test}, t)} < \eps
\end{equation}
then the split line $ t_{k}$ is identified as unnecessary and discarded from the set of split lines $t_{split} $.

For the second test, called the $\del$ test, the necessary length of a subinterval to retain a split line selection is a lower bound tolerance $\del$, which is heuristically chosen as a fraction of the data time interval.  As in the first test, an illustration of whether or not a split line is retained or discarded is portrayed by the following figures:
\begin{center}
\begin{minipage}{.5\textwidth}
\centering
\begin{tikzpicture}
 \draw[->] (4.0, -1.7) -- (8.0, -1.7) node[right] {$t$};
  \draw[->] (4, -1.7) -- (4, 1.3) node[above] {$z$};
\filldraw[black] (6,-1.7) circle (2pt) node[anchor=north] {$t_k$};
\filldraw[black] (5,-1.7) circle (2pt) node[anchor=north] {$t_{k-1}$};
\filldraw[black] (7,-1.7) circle (2pt) node[anchor=north] {$t_{k+1}$};
\draw (5.5,1.0) node[anchor=south] {$> \del$};
\draw[gray, thick] (5,-1) -- (6,-1) -- (6,1) -- (5,1) -- (5,-1);
\draw[gray, thick] (6,-1.4) -- (7,-1.4) -- (7,0.8) -- (6,0.8) -- (6,-1.4);
\end{tikzpicture}
	\captionof{figure}{Retained Split Line $t_k$}
 	\label{fig:splitdel1}
\end{minipage}%
\begin{minipage}{.5\textwidth}
\centering
\begin{tikzpicture}
 \draw[->] (4.0, -1.7) -- (8.0, -1.7) node[right] {$t$};
  \draw[->] (4, -1.7) -- (4, 1.3) node[above] {$z$};
\filldraw[black] (6,-1.7) circle (2pt) node[anchor=north] {$t_k$};
\filldraw[black] (5.5,-1.7) circle (2pt) node[anchor=north] {$t_{k-1}$};
\filldraw[black] (7,-1.7) circle (2pt) node[anchor=north] {$t_{k+1}$};
\draw (5.75,1.0) node[anchor=south] {$< \del$};
\draw[gray, thick] (5.5,-1) -- (6,-1) -- (6,1) -- (5.5,1) -- (5.5,-1);
\draw[gray, thick] (6,-1.4) -- (7,-1.4) -- (7,0.8) -- (6,0.8) -- (6,-1.4);
\end{tikzpicture}
	\captionof{figure}{Discarded Split Line $t_k$}
	 \label{fig:splitdel2}
\end{minipage}
\end{center}

\vspace{10pt} 

For an example of the $\del$ test, if the data range is $[0, 1]$, a chosen 10\% tolerance or $\del = 0.1$ subinterval length $[t_{k-1}, t_{k}]$ for $k \in  \set{1,2,\dotsc,{n-1}}$ is the minimum amount required to retain the split line. Otherwise, if the following relationship is satisfied:
\begin{equation}\label{del} \abs{t_{k}- t_{k-1}} < \del \end{equation}
then the split line $ t_{k}$ is identified as unnecessary and discarded from the set of split lines $t_{split} $. The value of the second subinterval length $\abs{t_{k+1}- t_{k}}$ stemming from the $\eps$ test is only evaluated for $k = n-1$, since each adjacent subinterval pair is evaluated across the entire data time interval proceeding from the initial time $t=0$ to the final time $t=T$.


Recursion is performed until a sufficient number of iterations is reached. In the case that all chosen splits $t_{k}$ for $k \in \set{1,2,\dotsc,n-1}$ are discarded, the method restarts with a different selection of $t_{split} $ values by adding in one more split than was used in the previous iteration. Consequently, for data displaying one distinct pattern across the time interval, it is expected that no split lines within the interval are necessary to be retained as long as the $\eps$ and $\del$ tests are satisfied to the given tolerances.

When a split line $t_{k}$ is retained, or the split line $t_{k}$ has ``passed" the $\eps$ and $\del$ tests, the algorithm is run recursively within the subinterval $[t_{k-1}, t_{k}]$ in order to find any further distinct patterns in the data subinterval.

After completing the maximum number of iterations, the algorithm outputs $t_{split} $ with split lines that satisfy the $\eps$ and $\del$ tests for a chosen point in space $x$. For added robustness, the algorithm is run for a number of different $x_{test}$ values, and the resulting $t_{split}$ with split lines that consistently satisfy the $\eps$ and $\del$ tests for these spatial points are chosen to define the subintervals for the DMD routine that follows.



A more explicit breakdown of the $n$-split algorithmic process is given by the following pseudo-code:
\begin{center}
\captionof{table}{$n$-split Algorithm}
\rule{\textwidth}{0.5pt}
 \end{center}
\begin{algorithmic}
    	\INPUT{dataset $Z$, split lines $t_{split}$, maximum iterations $M$}
    	\OUTPUT{split lines $t_{split}$, iteration count}
    	\WHILE{ iterations $< M$ }
    	\STATE{Choose $x_{k}$ randomly from the spatial nodes of $Z$}
    	\STATE{Set $s$ as number of split lines}
    	\STATE{Set $\eps$ as fraction of $\norm{Z(x_{k},[0,T])}_{\infty}$ and $\del$ as fraction of $[0,T]$}
    	\FOR{$k = 1$ to $s-1$}
		\IF {\begin{align*}
		\max \set{\bmm
		\abs{\max Z(x_{k},[t_{k-1}, t_{k}] ) - \max Z(x_{k},[t_{k} t_{k+1}] )}  \\
		\abs{\min_{t \in [t_{k-1}, t_{k}] } Z - \min_{t \in [t_{k}, t_{k+1}] }Z} 
		\emm } > \epsilon 
		\end{align*} }
			\IF{ $\abs{t_{k}- t_{k-1}} > \del$ }
				\STATE{Retain split line $t_k$ and add a new split line in $[t_{k-1}, t_{k}] $}
				\STATE{Run $n$-split algorithm on subinterval $Z(x,[t_{k-1}, t_{k}])$ and add new split lines as required}
				\STATE{Increase iteration count by output in previous step}
			\ELSE 
				\STATE{return}
    			\ENDIF
			\IF{ $k=s-1$ and $ \abs{t_{k+1}- t_{k}} > \del$}
				\STATE{similar steps as previous if statement}
    			\ENDIF
    		\ENDIF
     \ENDFOR
    	\IF{no split lines are retained} 
		\STATE{Increase iteration count by 1 and choose a new guess for $t_{split}$}
		\STATE{Run $n$-split algorithm on $Z $ and replace current split lines by output}
		\STATE{Increase iteration count by output in previous step}
    	\ENDIF
	\STATE{Discard redundant split lines}
    \ENDWHILE
\end{algorithmic}
\rule{\textwidth}{0.5pt}
\vspace{1pt} 

Thus, stepping out of the $n$-split algorithm, the split DMD algorithmic process is given by the following pseudo-code:

\begin{center}
\captionof{table}{Split DMD Algorithm}
\rule{\textwidth}{0.5pt}
 \end{center}
\begin{algorithmic}
    	\INPUT{rank of approximation $r$, dataset $Z$}
    	\OUTPUT{Order-$r$ approximations $Z_{DMD}$, Order-$r$ DMD modes $\Phi^{DMD}$}
    	\STATE{Choose split lines $t_{split}$ heuristically and set maximum number of iterations $M$ }
    	\FOR{enough iterations to reach steady-state solution}
    		\STATE{Run $n$-split algorithm on $Z$}
      \ENDFOR
	\STATE{Retain steady-state selection of split lines $t_{split}$}
    	\STATE{Set $s$ as number of split lines}
    	\FOR{$k = 1$ to $s$}
		\STATE{Run DMD algorithm to output DMD modes $\Phi_{k}^{DMD}$ 
		\\ and rank-$r$ $Z_{DMD}(x,[t_{k-1}, t_{k}] )$}
     \ENDFOR
		\STATE{Set $\Phi^{DMD} = \set{\Phi_{k}^{DMD}}_{k=1}^s$ }
		\STATE{return  $Z_{DMD},\Phi^{DMD}$}
\end{algorithmic}
\rule{\textwidth}{0.5pt}

\section{Results}

The main purpose of this study is to compare the effectiveness of the split DMD method, i.e. selecting DMD modes in subdomains of the entire solution, relative to DMD methods that select DMD modes across the entire domain on a nonlinear PDE problem. A simplified version of the proposed $n$-split algorithm in Section 2 was implemented for this study, by using an evenly-spaced initial guess for the split lines and a heuristic choice for the number of initial split lines.

\subsection{Comparison of Split OD and OD}

The split DMD algorithm was applied using both OD and standard DMD to the following form of the Kuramoto-Sivashinsky (KS) equation:
\begin{align}\label{KS1}
 \frac{\pde w}{\pde t} + \frac{\pde^4 w}{\pde x^4} = - 2w \frac{\pde w}{\pde x} -  \frac{\pde^2 w}{\pde x^2}
\end{align}
The stabilizing fourth-order and destabilizing second-order terms in the KS equation mimic the Navier-Stokes' energy behavior \cite{podbook}. An important feature of the KS equation is the \bol{bifurcation parameter} $L$, the length of the periodic domain where the model is studied, that determines the behavior of the dynamical system.
To perform a non-dimensionalization of the KS equation to a unit periodic domain $\bar x = \frac{x}{L}$, define $\eps = \frac{1}{L^2}$ so that the non-dimensional model is
\begin{align}\label{KS2}
  \frac{\pde w}{\pde t} +\eps^2 \frac{\pde^4 w}{\pde  \bar x^4} = - 2\eps w \frac{\pde w}{\pde \bar  x} - \eps \frac{\pde^2 w}{\pde \bar  x^2}
\end{align}
This exposes the relationships between length versus nonlinearity and stabilizing versus destabilizing terms.
Solutions to the non-dimensional equations above are considered in the domain
\[ (\bar x,t) =  \mathbb{T} \times (0,T) \]
with periodic domain $\mathbb{T}=(0,1)$, different values of the length parameter $L$ (and thus $\eps$) and the final time $T$.
These equations are solved with periodic boundary conditions
\begin{align*}
 w(0,t) = w(1,t),~~ w_{{\bar x}}(0,t) = w_{{\bar x}}(1,t),\\
 ~~
 w_{\bar x\bar x}(0,t) = w_{\bar x\bar x}(1,t),~~ w_{\bar x\bar x\bar x}(0,t) = w_{\bar x\bar x\bar x}(1,t)
\end{align*}
and initial condition
\begin{align}
w(\bar x,0) = \frac{\sin(4\pi \bar x)}{\sqrt{\eps}}
\end{align}
The dynamical behavior of a system modeled by the KS equation is summarized in the following table \cite{podbook}.

\begin{center}
\captionof{table}{Dynamical Behavior of Kuramoto-Sivashinsky}
\begin{tabular}{lll}
\hline
Length $L$ & \multicolumn{1}{c}{$\eps$} & Solution Behavior \\
\hline
12.5664 $\approx 4\pi$  &  0.00633257  & Bifurcation \\
12.8767  &  0.00603102  & Heteroclinic Bifurcation \\
13.1403  &  0.00579148  & Hopf Bifurcation \\
402.2590  $\approx 128\pi$ &  0.00000618  & ``Chaos'' \\
\hline \\
\end{tabular}
\end{center}

The periodic domain lengths $L$ initially tested were 12.6, 13.2, and 402.3.
The following figures portray a visual inspection of the improvement of reduced-order approximation to the full model due to the splitting procedure of the split DMD algorithm:

\begin{center}
\begin{minipage}{.5\textwidth}
  \centering
\includegraphics[scale=0.4]{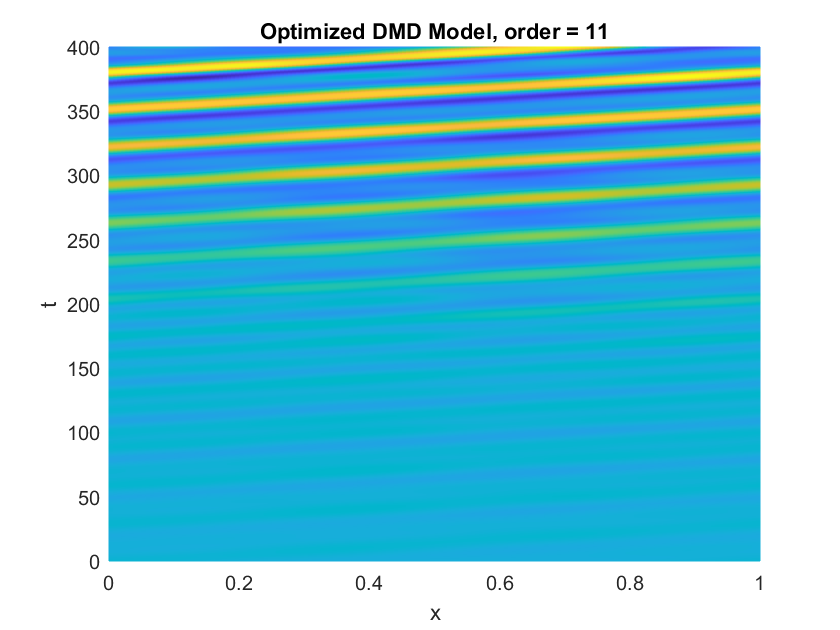}
	\captionof{figure}{KS OD, $L = 12.6$}
 	\label{fig:od126}
\end{minipage}%
\begin{minipage}{.5\textwidth}
  \centering
\includegraphics[scale=0.4]{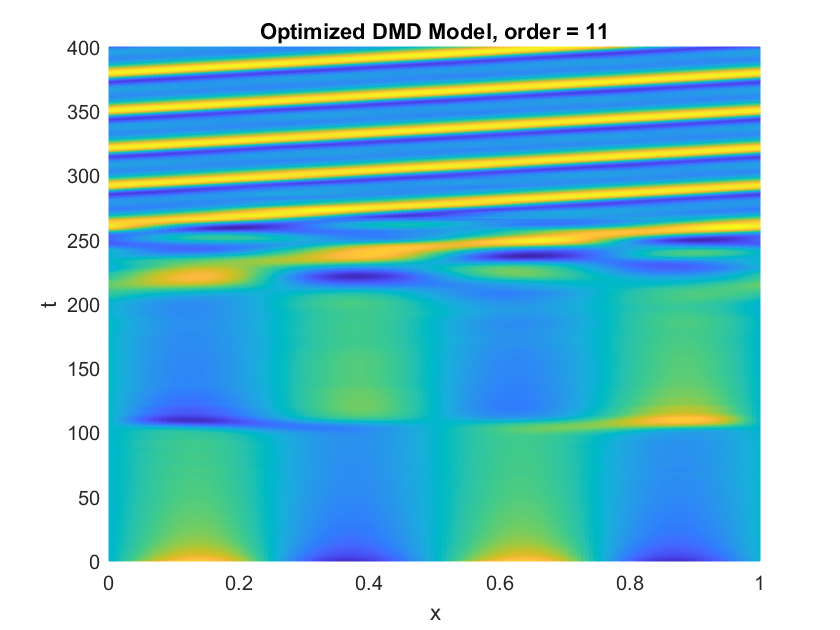}
	\captionof{figure}{KS 4-split OD, $L = 12.6$}
	\label{fig:sod126}
\end{minipage}
\end{center}

\begin{center}
\begin{minipage}{.5\textwidth}
  \centering
\includegraphics[scale=0.4]{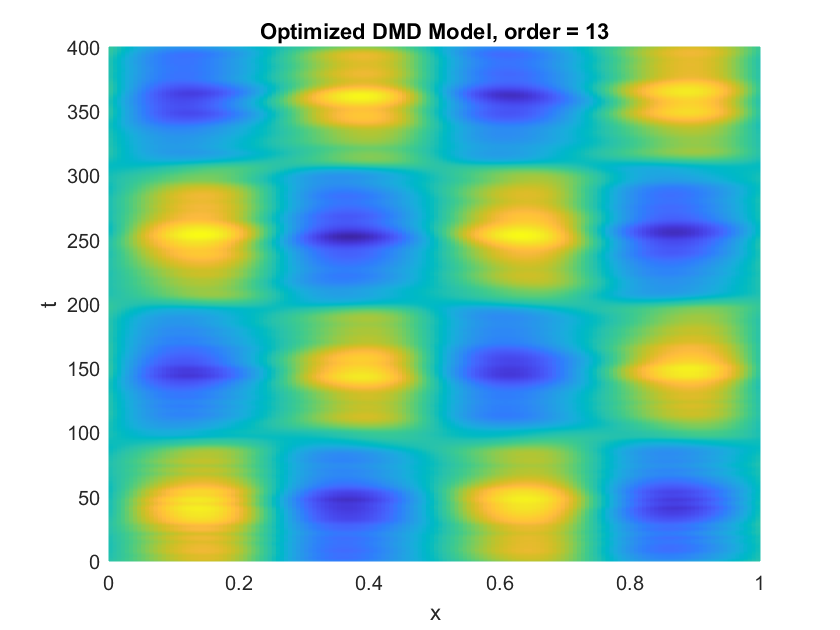}
	\captionof{figure}{KS OD, $L = 13.2$}
 	\label{fig:od132}
\end{minipage}%
\begin{minipage}{.5\textwidth}
  \centering
\includegraphics[scale=0.4]{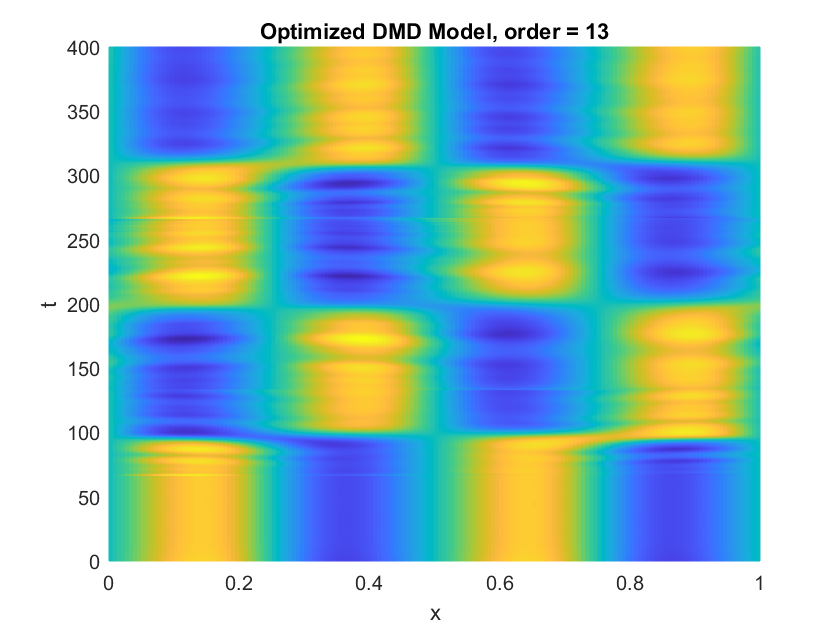}
	\captionof{figure}{KS 4-split OD, $L = 13.2$}
	\label{fig:sod132}
\end{minipage}
\end{center}

\begin{center}
\begin{minipage}{.5\textwidth}
  \centering
\includegraphics[scale=0.4]{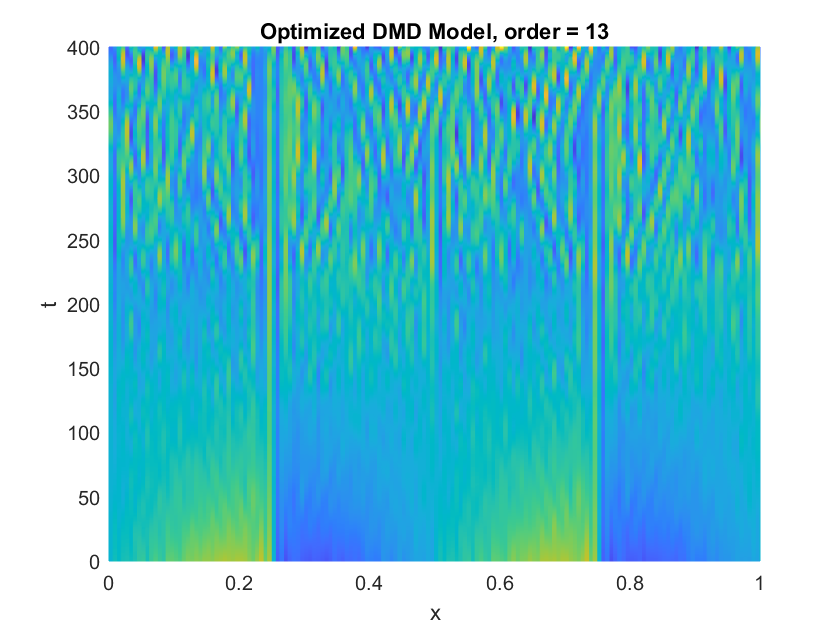}
	\captionof{figure}{KS OD, $L = 402.3$}
 	\label{fig:od4023}
\end{minipage}%
\begin{minipage}{.5\textwidth}
  \centering
\includegraphics[scale=0.4]{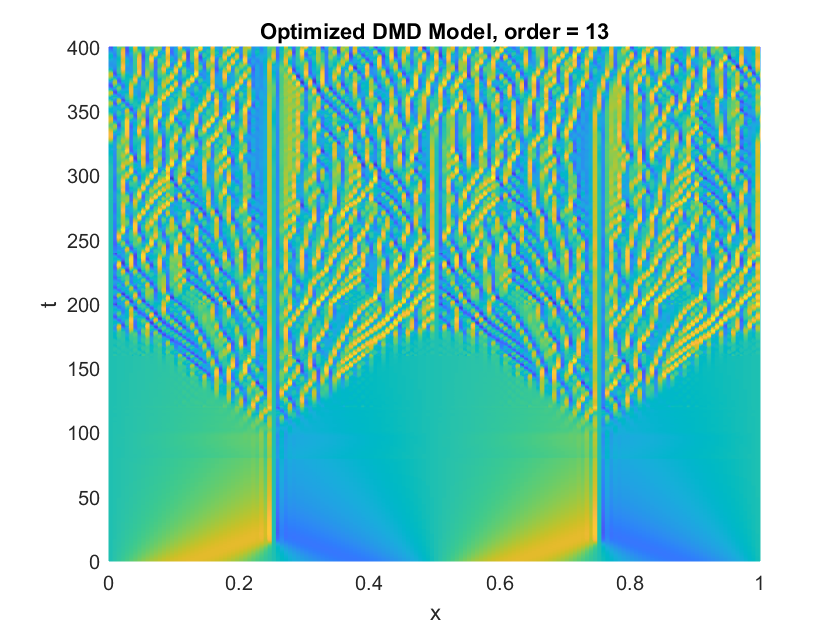}
	\captionof{figure}{KS 10-split OD, $L = 402.3$}
	\label{fig:sod4023}
\end{minipage}
\end{center}

The relative error of the split OD models provided in the preceding figures shows a clear improvement in comparison to the standard OD models with respect to the finite element solution. Without splitting, it is evident that OD does not provide a representation of the data that captures the coherent structures throughout the time interval.

The following table provides a comparison of the split OD ROMs versus the OD ROMs for $t=400$, with information on each reduced-order model concerning:
\been
\item Number of splits used: $n$
\item Length of the periodic time interval (bifurcation parameter): $L$
\item Rank of the approximation: $r$
\item Computation error: $\norm{r(t)}_2$
\item Computation time (in seconds)
\enen

\begin{center}
\captionof{table}{Split OD ROM vs OD ROM, $t=400$}
\begin{tabular}{c l c l c l c l c l c}
\hline
 $n$-Split & Length & Rank & Computation error  & Time \\
\hline
0   & 12.60 & 11 & 0.539273074481482 & $1.656649834 \times 10^2$ \\
4  & 12.60 & 11 & 0.086861903201658 & $1.424075578 \times 10^2$ \\
\hline
0  & 13.20 & 13 & 0.450356679589967 & $2.765350954 \times 10^2$ \\
4   & 13.20 & 13 & 0.187951801209479 & $1.874434272 \times 10^2$ \\
\hline
0  & 402.3 & 13 & 18.071816203629595 & $3.078671183 \times 10^2$ \\
10   & 402.3 & 13 & 2.266344147532708 & $1.270328812 \times 10^2$  \\
\hline \\
\end{tabular}
\end{center}

The table suggests a substantial improvement in using the splitting approach.

\subsection{Sensitivity of Initial Condition}

For the chaotic dynamics ($L=402.3$), the sensitivity of the results in Section 4.1 was tested by adding a random perturbation to the discretization of the previous initial condition. Letting $\beta(x) = \text{unif}\paren{0,\frac{1}{20}}$, a uniformly distributed random number at each spatial location (which is 161 in this case), the new initial conditions are:
\begin{align}\label{ICKSp}
w(x,0) = \frac{\sin(4\pi x)}{\sqrt{\eps}} + \beta(x),~~w_x(x,0) = \frac{4\pi\cos(4\pi x)}{\sqrt{\eps}}
\end{align}

The results 
are summarized in the following table, compared to the standard OD and DMD models 
:
\begin{center}
\captionof{table}{Split DMD ROM vs Split OD ROM, $\beta$ shift}
\begin{tabular}{c l c l c l c l c l c | c}
\hline
 ROM & $n$-Split & Length & Rank & Computation error  & Time \\
\hline
DMD & 0  & 402.3 & 13 & 28.927186508313124 & 0.1634884 \\
DMD & 10  & 402.3 & 13 & 5.798752827124941 & 0.17863 \\
\hline
OD & 0   & 402.3 & 13 & 22.233228342958030 & $2.616457798 \times 10^2$  \\
OD & 10   & 402.3 & 13 & 2.739207339916917 & $1.858652532 \times 10^2$  \\
\hline
\end{tabular}
\end{center}


\subsection{Sensitivity of Bifurcation Parameter}

For the chaotic dynamics ($L=402.3$), the sensitivity of the results in Section 4.1 was tested by shifting the bifurcation parameter to $L=402.35$.

The results 
are summarized in the following table:
\begin{center}
\captionof{table}{OD ROM vs Split OD ROM, $L=402.35$}
\begin{tabular}{c l c l c l c l c l c }
\hline
 $n$-Split & Length & Rank & Computation error  & Time \\
\hline
0   & 402.35 & 13 & 19.038661517042843 & $2.773977651 \times 10^2$  \\
10   & 402.35 & 13 & 2.670265429612690 & $1.323419021 \times 10^2$  \\
\hline
\end{tabular}
\end{center}


\subsection{Sensitivity of Split Lines}

For the chaotic case in Section 4.1, the sensitivity of the results was tested by shifting the split lines. The four sensitivity test cases for the 10-split OD model are uniform shifts by $\pm 1, \pm 3$ seconds in each of the interior time splits.

The results 
are summarized in the following table, with the ``Shift" column values measured in seconds:
\begin{center}
\captionof{table}{Split OD ROM with split line shifts}
\begin{tabular}{c l c l c l c l c l c l c }
\hline
 $n$-Split & Shift & Length & Rank & Computation error  & Time \\
\hline
10 & +3 & 402.3 & 13 & 3.362777585745632 & $1.464054566\times 10^2$  \\
10  & +1 & 402.3 & 13 & 4.197440908311128 & $1.340162133 \times 10^2$  \\
10  & -1 & 402.3 & 13 & 5.702330507412169 & $1.3498605 \times 10^2$  \\
10  & -3 & 402.3 & 13 & 4.645326326310565 & $1.42945084 \times 10^2$  \\
\hline
\end{tabular}
\end{center}


\section{Discussion}

As mentioned in Section 2, the OD method without splitting has issues reconstructing systems with bifurcation present. 
Using the $n$-split algorithm, 
the split OD method reconstructs the system in question to a level of accuracy sufficient for recognizing the effects of bifurcation on the system dynamics. 

The standard DMD algorithm does not provide any useful results in the test case in Section 4.1.
By the $\beta$ shift in Section 4.2, it is evident that the split DMD is able to provide results that sufficiently reconstruct the system dynamics in certain cases. The anticipation 
is that the split DMD, although computationally cheap, is limited in its potential solution accuracy. In any case, accuracy of a DMD solution which rivals an OD solution for the KS equation is not expected.

An increase in periodic length affects the solution to the standard OD ROM significantly, whereas the 10-split OD model maintains a high level of accuracy. Furthermore, the robustness to where the split occurs in the chaotic dynamics 
implies that a sufficient number of splits is enough to produce an accurate reconstruction of the finite element solution.

The main issue with the splitting approach is the inability to predict future time states with the resulting modes. 
In the $L=13.2$ case, the solution for each split can simply be ``copied" in anticipation of the evident recurring pattern. However, in the $L=402.3$ case, the chaotic dynamics (implying differing coherent structures for each split) do not provide a gateway to predict the the chaotic pattern.

The conclusion of these results is that a decrease in computation time and increase in accuracy of solutions suggests that the $n$-split algorithm for DMD methods is superior to the standard DMD and optimized DMD methods for reconstructing a solution to the KS equation.

The numerical study applying the split DMD method to the KS equation suggests that there is a benefit to using this splitting approach for modeling dynamical systems with certain structural features. However, the work described in this text did not factor in many possibilities that may come from studying dynamical systems theory. Papers such as \cite{kass02} provides a method to choose modes over a longer time period that are not biased by the transient dynamics. Page 52 of \cite{dmdbook} provides a numerical method for inspecting the attracting manifolds of the dynamical system. \cite{mezicbook} and Chapter 3 of \cite{dmdbook} provide direction towards future state prediction based on Koopman theory. Further study on these topics can lead to a more sophisticated approach towards creating a metric for attractors when chaotic dynamics are involved.

\end{document}